\documentclass[letterpaper, 10 pt, conference]{ieeeconf}
\IEEEoverridecommandlockouts
\usepackage{cite}
\usepackage{nicefrac}

\usepackage{ulem}
\usepackage{xcolor}
\usepackage{multirow}
\usepackage{enumerate}
\usepackage{booktabs} 
\usepackage{amsmath}
\usepackage{physics}
\usepackage{hyperref}
\usepackage{url}
\usepackage{graphicx}
\usepackage{siunitx}
\usepackage{import}
\usepackage{placeins}
\usepackage{amsthm}
\usepackage{amssymb}
\usepackage{empheq}
\usepackage{ulem}
\usepackage{algorithm}
\usepackage{float}
\usepackage{multicol}
\usepackage{cleveref}
\usepackage[english]{babel}

\newtheorem{theorem}{Theorem} 
\theoremstyle{definition}

\newtheorem{lemma}{Lemma}

\newcommand{\hide}[1]{}

\newcommand{\eop}{\hfill{$\square$}}

\newcommand{\Om}{\Phi}
\newcommand{\om}{\phi}
\newcommand{\ups}{\Upsilon}
 \newcommand{\Ups}{{\mathbf \Upsilon}}

\newcommand{\offs}{\Gamma}

\newcommand{\Cx}{C}
\newcommand{\cx}{c}

\newcommand{\bax}{{\bar a}}
\newcommand{\Ax}{A}
\newcommand{\ax}{a}

\newcommand{\pc}{\psi^c}
\newcommand{\pa}{\psi^a}

\newcommand{\Pc}{\Psi^c}
\newcommand{\Pa}{\Psi^a}
\newcommand{\bm}{\bar{m}}

\newcommand{\cd}{{\cal D}}

\newcommand{\arxiv}[2]{#2}  

\usepackage{empheq}
\newcommand*\widefbox[1]{\fbox{\hspace{0em}#1\hspace{0em}}}

\newcommand{\floor}[1]{\lfloor #1 \rfloor}

\usepackage{textcomp}
\def\BibTeX{{\rm B\kern-.05em{\sc i\kern-.025em b}\kern-.08em
    T\kern-.1667em\lower.7ex\hbox{E}\kern-.125emX}}

\begin{document}

\title{Saturated total-population dependent branching process and viral markets
}

\author{Khushboo Agarwal* and Veeraruna Kavitha\\
\textit{\{agarwal.khushboo, vkavitha\}@iitb.ac.in}, IEOR, IIT Bombay,  India
\thanks{*The work of first author is partially supported by Prime Minister's Research Fellowship (PMRF), India.}
}

\maketitle

\begin{abstract}
Interesting posts are continually forwarded by the users of the online social network (OSN). Such propagation leads to re-forwarding of the post to some of the previous recipients, which increases as the post reaches a large number of users. Consequently, the effective forwards (after deleting the re-forwards) reduce, eventually leading to the saturation of the total number of copies. We model this process as a new variant of the branching process, the `saturated total-population-dependent branching process', and analyse it using the stochastic approximation technique. Notably, we obtain deterministic trajectories which approximate the total and unread copies of the post `asymptotically and almost surely' over any finite time window; this trajectory depends only on four parameters related to the network characteristics. Further, we provide expressions for the peak unread copies, maximum outreach and the life span of the post. We observe known exponential growth but with time-varying rates. We also validate our theory through detailed simulations on the SNAP Twitter dataset.
\end{abstract}

\section{Introduction}
Social media offers a global platform for communication;  people share their content, forward interesting posts among the ones shared with them, etc. It plays a significant role in marketing, e.g., a post advertising a specific product/service can get viral (reach a large number of users). It is also used for spreading propaganda, sharing a piece of knowledge, influencing the beliefs of a large population, etc.

The content propagation (CP) over an online social network (OSN) can be outlined as follows: a) a post is shared with an initial set of users, called seed users; b) the seed users forward the post to their friends/followers if they like it;  c) the recipients follow suit and this continues. The content either gets extinct in the initial phase or gets viral, and eventually, copies get \textit{saturated}, and the propagation ceases.

There are several approaches for studying CP, and we use branching processes (BPs) to analyze the same. Further, since a post can witness a huge surge in the shares in a short duration, we consider continuous-time BPs (e.g., \cite{dhounchak2017viral, iribarren2011branching, van2010viral}). We discuss and compare other techniques at the end. 

\noindent \textbf{Branching processes:} A variety of BPs have been studied before (e.g., \cite{athreya2001branching,agarwal2021co,klebaner1984population} are few strands of  them). We discuss a few relevant varieties and provide some details to model the CP process using BPs. 
 The course of any BP largely depends upon the expected number of offsprings ($\Gamma$).
 We have super-critical BPs   when the  
expected offsprings $E[\Gamma] > 1$ (see e.g., \cite{athreya2001branching}); such BPs have positive probability of exploding (population  grows exponentially) and  can mimic the viral CP.    In the critical/sub-critical regime ($E[\Gamma]\leq1$), the BP gets extinct with probability one, i.e., the population eventually declines to zero. The current population in a BP represents the live population, while the total population also includes the dead ones. In the context of  CP, the number of unread copies (the corresponding recipients are yet to view/read the post) represents the current population. While the total number of recipients, including the ones that already read the post, represents the total population.

\noindent \textbf{Saturation:} When a post is viral and is already forwarded to a noticeable fraction of the network, a significant fraction of further forwards (by future users) can overlap a  part of the network that already received the post. The expected number of effective forwards (after deleting re-forwards) represents the expected offsprings when one attempts to model CP using a BP. The number of re-forwards depends on total copies (number of users that already received the post) and not just on the current population/unread copies.
We are only aware of current-population  dependent BPs (e.g., \cite{klebaner1984population}).

Thus the existing BP models are insufficient to mimic a saturated CP due to two imperative factors: a) the overlaps in forwards/offsprings depend upon the total copies/population, and b) the usual content propagation process traverses from super-critical to sub-critical regime before getting extinct. 
We thus consider a  new variant of continuous-time total-population dependent Markovian BP, named as \textit{saturated total-population dependent BP (STP-BP)}. The total copies either increase with time or saturate; if re-forwards are proportional to the total copies, then the resultant value of expected offsprings decreases with time. This also ensures the desired transition from super-critical to sub-critical regime. The saturated BP well mimics the CP.

Recently, we analyzed many new variants of BPs, using a new approach based on stochastic approximation techniques, e.g., attack and acquisition  BP (competing viral markets \cite{agarwal2021co}), and proportion-dependent BP (fake news on OSNs \cite{kapsikar2020controlling}). Further modifying the said approach to address the required finite horizon analysis allows us to analyze saturation resulting from total population dependency.

\noindent \textbf{Key contributions:} At first, we formally analyze the STP-BP. \textit{We derive an appropriate ordinary differential equation (ODE) and its solution, which (time asymptotically) almost surely approximates the embedded chain of the CP process over any finite time window. } These deterministic solutions depict exponential growth and linear fall for unread copies. Secondly,  we model and fit an appropriate total-copies dependent piece-wise linearly decreasing function for the expected offsprings, having two different slopes. 
We further derive important metrics like the peak number of unread copies, time asymptotic value of the total copies, and others. The growth of the total copies is exponential and depends on the reduction rates in the expected forwards.

We corroborate our theoretical results by performing Monte-Carlo simulation on SNAP Twitter-dataset \cite{mcauley2012learning}. 
The description of the theoretical trajectories depends only on four parameters of the OSN (e.g.,  two rates of reduction of the expected forwards) and the attractiveness of the post.
  
\noindent \textbf{Related work:} The study of CP on OSNs has been a topic of interest for a long time, and several approaches are used for its analysis. 
Random graph models are widely used to analyze CP on OSNs (e.g., see \cite{deijfen2016winner, zhou2019cost, lu2019beyond}); in particular, \cite{lu2019beyond} considers re-exposure of users with a post on OSN. However, such models can not capture aspects like virality, as discussed in \cite{van2010viral,dhounchak2017viral, iribarren2011branching}, while  BPs facilitate virality analysis like growth patterns, extinction/virality probability etc.

The epidemiology-based models (in particular SIR) are also used to study CP, which succeed in capturing saturation (e.g., \cite{rodrigues2016can, freeman2014viral, jiang2019quantitative}). The current and total copies are respectively modelled as infected and (infected$+$recovered) populations. We argue that this approach is not suitable for analyzing viral markets on OSNs/email platforms: (a) this set of literature directly starts with an appropriate ODE (e.g.,  \cite{freeman2014viral} considers exponentially diminishing infection rate/interest in the post, while \cite{rodrigues2016can} considers standard SIR ODE), (b) they do not delve into the details of the random dynamics (e.g.,  the chance encounters between various individuals), and more importantly, (c) on OSNs,  majority of users share the post to a subset of their friends  `only once' after viewing and lose interest immediately after; in contrast, in SIR based models, infected individuals keep infecting/spreading the information for a random/prefixed time before recovering/losing-interest. In other words, \textit{SIR-based models well capture the typical behaviour of word-of-mouth dynamics (individuals remain interested in gossip,  keep sharing it, and then lose interest) and are insufficient for viral markets over OSNs}.

There is a brief indirect mention of the saturation effect in \cite{van2010viral}, where using BPs, the authors predict the future progress of CP using the available history of an ongoing campaign; they use the well-known Kolmogorov's backward equations  (for PGFs in BPs) to achieve this. In contrast, we provide a theoretical study of a new relevant variant of BP, facilitating an exhaustive study of the saturated CP on OSNs.

\vspace{-1mm}
\section{Problem description and background}\label{sec_model}
Consider an OSN where the content of interest is forwarded by its users. At the start of the propagation, the post is shared by its content provider to  an initial set of users,  called seed users.
These users view the post on  the OSN at random time instances. Then, they forward it to some or all of their friends, depending on how much they like the post.\footnote{The post may seem appealing to the users for the offers mentioned, the creativity or the informational quotient of the content (see \cite{van2010viral}).} This subset of users further forward the post when they visit the OSN\footnote{When a user views, reads and forwards some of the posts on its timeline.}, and the post propagation continues likewise (more details in  \cite{dhounchak2017viral, agarwal2021co, van2010viral}). 
The random instances at which users visit the OSN are called wake-up times. The time gap between the wake-up and the post-reception times of any recipient  is  exponentially distributed with parameter $\lambda$.  
Let $\mathcal{F}$ be the (random) number of friends of a typical user, with finite mean. 
The attractiveness of the post  is captured by  factor $\rho$, which specifies a subset of friends $(\offs)$ to whom the post is forwarded. Naturally, $\rho$ determines the  growth of the post.

We will see that the content propagation (CP) majorly depends on the expected number of such forwards, which in turn depends on the then total recipients of the post (say $a$). That is, $E[\offs] := m(a)$, where $m$ is some function which will be discussed in detail in the coming sections. The study of such a dependency (using SNAP dataset) and its influence on the CP is the key differentiating feature of our work. As one may anticipate again, such dependency will eventually lead to the reduction in the effective forwards of the post, in turn leading to the saturation of  shares. This effect has been  observed majorly through numerical studies in the past, e.g.,  see \cite[Fig. 6]{dhounchak2017viral} and \cite[Fig. 1]{hemsley2016studying} for total and current shares, respectively. \textit{We aim to provide an analytical explanation of saturation in terms of an approximating deterministic trajectory and some relevant performance measures.} We next describe the new variant of `saturated' BP, STP-BP.


\vspace{-2mm}
\subsection{Saturated Content propagation and Branching process}

Let $\Ax(t)$ denote the   total number of copies of the post on the network, i.e., all the copies which has been received (viewed or not, forwarded or not) by the users on the OSN till time $t$; we briefly refer this number as \textit{total shares}. Further, let $C(t)$ be the unread/live/current number of copies of the post till time $t$. These represent the copies that have been received but not yet viewed  by the users till time $t$, and are solely responsible for further propagation of the post; we refer to this number as \textit{current shares}.

Let $\Ax(0) = \Cx(0) = \ax_0$ be the number of seed users. By previous discussion, it is clear that the number of forwards have a direct dependence on the total shares, i.e., $\offs = \offs(A(t))$. Observe that $\offs$  does not depend on the current shares, as is usually considered in the BP literature.

Let $\tau^+$ and $\tau^-$ denote the time instances immediately after and before wake-up time $\tau$ (of any user),  e.g., $\Ax(\tau^+) := \lim_{t \downarrow \tau} \Ax(t)$. 
Then, the evolution of the system at transition epoch $\tau$ is:

\vspace{-5mm}
{\small\begin{equation}\label{evolve_cont}
\begin{aligned}
\hspace{1.5cm}\Cx(\tau^+) &= \Cx(\tau^-)  + \offs(A(\tau^-)) - 1, \\  
\hspace{1.5cm}\Ax(\tau^+) &= \Ax(\tau^-)  + \offs(A(\tau^-)).
\end{aligned}
\end{equation}}

Above dynamics can easily be placed  in a BP framework, when the current shares are modelled as a population. Accordingly, the number of forwards $(\offs(A(\tau^-)))$ can be viewed as the offsprings of the population. Further, when the post is shared by a user, the current shares reduce by $1$ (see \eqref{evolve_cont}); this is  exactly like a death in a BP (see \cite{dhounchak2017viral, kapsikar2020controlling, agarwal2021co} for similar details). 

To analyse total population dependent BPs, one needs to study two-dimensional tuple $\Om(t) := (\Cx(t), \Ax(t))$ simultaneously, a realisation of which is denoted by $\om = (\cx, \ax)$. In contrast, the existing BP models can analyse  $\Cx(t)$ alone if the need is only to analyse the current population. We now proceed towards analytically deriving the trajectories of the tuple $\Om(t)$, and other salient features.

\vspace{-1mm}
\section{Dynamics and ODE approximation }\label{sec_dynamics}

To facilitate the study of STP-BP, we analyse the embedded chain corresponding to the underlying continuous-time jump process (CTJP), as in \cite{kapsikar2020controlling, agarwal2021co}. It is a standard technique to use embedded chains, when transience, recurrence, extinction and similar properties of CTJP are studied. We use it for similar purpose in this paper. 
In particular, we observe the  dynamics in \eqref{evolve_cont} at the time instances when a user with unread copy of the post wakes-up. 
Let $\tau_n$ be $n^{th}$ such transition epoch.\footnote{If the post gets extinct at $n^{th}$ epoch, we set $\tau_{k} :=\tau_n$ for all $k \geq n$ and the same is true for rest of the quantities.}
Let $\Cx_n := \Cx(\tau_n^+)$ be the current shares of the post immediately after $\tau_n$. Similarly  define $\Ax_n$. Note that the  time taken by the first user to wake-up after $n^{th}$ transition epoch, ($\tau_{n+1}-\tau_n$), is exponentially distributed with parameter  $\lambda \Cx_n$. Thus, if a user wakes up at $\tau_n$, then:

\vspace{-3mm}
{\small
\begin{align}\label{evolve_SA}
\Cx_{n} &= \Cx_{n-1}  + \offs_{n}(A_{n-1}) - 1, \ \Ax_{n} = \Ax_{n-1}  + \offs_{n}(A_{n-1}).
\end{align}}

As in \cite{kapsikar2020controlling, agarwal2021co}, we use stochastic approximation (SA) approach to study the embedded chain. Towards this, define the following fractions of current and total shares respectively, $\Pc_n := \nicefrac{\Cx_n}{n}$ and $\Pa_n := \nicefrac{\Ax_n}{n}$ for $n \geq 1$, with $\Pc_0 = \Pa_0 := \ax_0$. Also let $\Ups_n := (\Pc_n,  \Pa_n)$. Further, define:

\vspace{-3mm}
{\small 
\begin{eqnarray}\label{eqn_eps_n}
    \epsilon_n = \frac{1}{n+1}, \ t_n := \sum_{k=1}^n \epsilon_{k-1}, \mbox{ and } \eta(t) := \max\left  \{ n: t_n \le t \right \}. \hspace{-2mm}
\end{eqnarray}}Then, the evolution of $\ups_n$ can be captured by 2-dimensional SA based updates given below (see \eqref{evolve_SA}):

\vspace{-0.25cm}
{\small
\begin{equation}\label{eqn_SA}
\begin{aligned}
\Pc_n &= \Pc_{n-1} + \epsilon_{n-1} \left[\offs_{n}(A_{n-1}) - 1 - \Pc_{n-1} \right ]1_{\Pc_{n-1} > 0}, \\
\Pa_n &= \Pa_{n-1} + \epsilon_{n-1} \left [\offs_{ n}(A_{n-1})  - \Pa_{n-1} \right ]1_{\Pc_{n-1} > 0}.
\end{aligned}
\end{equation}}

We analyse these fractions using SA techniques (e.g., \cite{kushner2003stochastic}), which helps in approximating the same using the solutions of the ODE (see \eqref{eqn_eps_n}):

\vspace{-3mm}
{\small
\begin{align}\label{eqn_ODE}
&\dot{\pc} = \left(m(\ax) - 1 - \pc\right) I, \  
\dot{\pa} = \left(m(\ax)  - \pa\right)I, \mbox{ with}\\
&I := 1_{\pc > 0}, \ \ax(t) :=     \pa(t) \eta(t), \mbox{ and } m(a) := E[\offs(\ax)]. \nonumber
\end{align}}

By Lemma \ref{lemma_existence}, the solution for this non-autonomous and non-smooth ODE exists over any finite time interval in the extended sense (satisfies the ODE for almost all $t$). In Theorem \ref{thrm1} given below,  we will prove that the above ODE indeed approximates \eqref{evolve_SA}.

\noindent \textbf{Approximation result:} The study of continuous time population size dependent BPs has been limited in the literature. 
In this paper, we use the ODE approximation result to study the saturated BP. Now, for mathematical tractability, we require the following assumption on offspring distributions: 
\begin{enumerate}
    \item[{\bf (A)}] There exists an integrable random variable, $\hat{\offs}$, such that $\offs(\ax) \leq \hat{\offs}$ almost surely for every $\ax$ and $E[\hat{\offs}]^2 < \infty$. 
    \item[{\bf (B)}] The mean function, $m(\cdot)$ is Lipschitz continuous.
\end{enumerate}The assumption \textbf{(A)} is readily satisfied (details in the next section), while \textbf{(B)} is an extra assumption required for additional affirmation (see Theorem \ref{thrm1}(ii)).
We will now see that the piece-wise constant interpolation, $\Ups^n(\cdot) := (\Psi^{n, c}(\cdot), \Psi^{n, a}(\cdot))$ of $\Ups_n$ trajectory defined as:
\begin{align}
\Ups^n(t) = \Ups_n \mbox{ if } t \in [t_n, t_{n+1}),
\end{align}
satisfies an almost integral representation  as below, with $I_n := 1_{\Psi^{n, c}(s) > 0}$ (see \eqref{eqn_linear2} for derivation):

\vspace{-4mm}
{\footnotesize
\begin{align}\label{eqn_piecewise_main}
\Psi^{n, a}(t) &= \Pa_n + \int_0^t \hspace{-1mm} \bigg(m(\Psi^{n, a}(s) n) - \Psi^{n, a}(s)\bigg)I_n ds + \varepsilon^{n, a}(t),\\
\Psi^{n, c}(t) &= \Pc_n + \int_0^t \bigg(m(\Psi^{n, a}(s) n) - 1 -  \Psi^{n, c}(s)\bigg)I_n ds  + \varepsilon^{n, c}(t). \nonumber
\end{align}}Let $\hat{\ups}^n(\cdot)$ be the  solution of ODE \eqref{eqn_ODE}, with $\hat{\ups}^{n}(0) = \ups_{n}$. Observe that $\Ups^n(\cdot)$ in \eqref{eqn_piecewise_main} is similar to $\hat{\ups}^n(\cdot)$, except for the difference term $\varepsilon^n(t) := (\varepsilon^{n, c}(t), \varepsilon^{n, a}(t))$. 
Further, if at all $\varepsilon^n(t) = 0$ for all $t \leq T $, then, by uniqueness of the  solution (see  Lemma \ref{lemma_existence}), the BP trajectory \eqref{eqn_SA} would have coincided with it, i.e., $\ups_k = \hat{\ups}(t_k)$ for all $k$ such that $t_k \leq T$. However, it is not true in general; nevertheless, we will show that $||\varepsilon^n|| \to 0$ as $n \to \infty$ (see norm $||\cdot||$ in \eqref{eqn_norm}). 

Thus, we have two operators which are converging towards each other; the first operator including $\varepsilon^n$ in \eqref{eqn_piecewise_main} provides the BP trajectory, while the second operator without $\varepsilon^n$ in \eqref{eqn_piecewise_main} provides the ODE solution.
Further, using Maximum theorem, we show that  the difference between  the two solutions of the operators \eqref{eqn_piecewise_main} with  and without $\varepsilon^n$ is small, when $||\varepsilon^n||$ is small. Formally, we state the result as follows:

\begin{theorem}\label{thrm1}
For any $\ups(\cdot) = (\pc(\cdot), \pa(\cdot))$, define the norm with any finite $T > 0$:
\begin{align}
\begin{aligned}\label{eqn_norm}
||\ups|| &:= \max\{||\pc||, ||\pa||\}, \mbox{ where } \\
||\psi^i|| &:= \sup\{t \in [0,T]: |\psi^i(t)|\} \mbox{ for any } i \in \{a, c\}.
\end{aligned}
\end{align}
Under assumption \textbf{(A)}, we have the following almost surely:
\begin{enumerate}
    \item[(i)] $||\varepsilon^n|| \to 0$ as $n \to \infty$, and 
    \item[(ii)] if \textbf{(B)} also holds, then, the difference  $\sup_{k : k \geq n, t_k \leq T} || \Ups^n(t_k) - \hat{\ups}^n(t_k) ||$ depends upon the magnitude of $||\varepsilon^n||$.
\end{enumerate}
\end{theorem}
\noindent \textbf{Proof} is provided in Appendix. \eop

\noindent \textbf{Remarks:} (i) For each $n$, consider the ODE initialised with value of the embedded chain, $\ups_{n}$. Then, the embedded chain values at transition epochs, $k \in [n , \eta(t_{n} + T)]$, are  close to the ODE solution, $\hat{\ups}^{n}(t_k-t_{n})$,  at time points $t_k \in [t_{n}, t_{n}+T] $. This approximation improves as $n$ increases, and \textit{the result is true almost surely} (a.s.) and for all $T < \infty$.

\noindent (ii) \textbf{Dichotomy:} We have a `modified 
dichotomy': either the population gets extinct ($\ups_{n} = 0$ in initial epochs), or the population explodes exponentially  as confirmed by  ODE-solution \eqref{eqn_total_ode}   in section \ref{ode_based_analysis}.  In contrast to classical dichotomy (e.g., \cite{athreya2001branching}), \textit{in both the sets of sample paths, the saturated BP eventually gets extinct}.

\noindent (iii) As  the ODE solution approximates the embedded chain (BP) trajectory, one can analyse the latter using the former.  \textit{In contrast to many existing studies, we would consider the analysis of the ODE trajectories and not the attractors, which is more relevant here. } Prior to that we derive an appropriate function that can represent $m(\cdot)$, the total-share-dependent expected forwards (TeF) in a typical OSN.

\vspace{-1mm}
\section{Population dependent expected forwards}\label{sec_mean}
\vspace{-1mm}

It is clear from the ODE \eqref{eqn_ODE} and Theorem \ref{thrm1} that the expected number of forwards, $m(\cdot)$, influence the course of any post; we construct a piece-wise linear function to capture the same. We estimate the parameters and validate the model using the SNAP dataset \cite{mcauley2012learning} in section \ref{sec_numerical}. Following aspects are considered for the model:

$\bullet$ Any user on the OSN forwards the post to a subset of friends based on its attractiveness.
Among this subset, a fraction of users will have previously received the post. 
\textit{Most likely, such users will not be interested in the same post again.} Say $\kappa$ is the fraction of common friends between any two typical users, then $\kappa \ax$ denotes the number of such re-forwards if `$\ax$' number of users already had the post. \textit{Thus a linearly decreasing function is a suitable choice for $m(\cdot)$.   }

$\bullet$ It is likely that a user with a high number of friends is also a friend to a large number of users; the dataset supports this observation. 
Thus, a user with more connections is more likely to receive the post earlier. 
To capture such aspects, {\it we model the friends of a user, ${\cal F}(a)$, to be independent across users, but with decreasing expected values, i.e., $a \mapsto E[{\cal F}(a)]$ is itself a decreasing function}. With this, {\it the assumption \textbf{(A)} is readily satisfied by considering for example, $\hat{\offs} = {\cal F}(\ax_0)$}, with $E[{\cal F}(\ax_0)]^2 < \infty$.

$\bullet$ Further, as time
progresses, the post would naturally proceed towards saturation. Therefore, it is likely that the slope of $m(\cdot)$ is drastically different towards the end. To summarize,
\textit{we model the TeF function $m(\cdot)$ as a piece-wise linearly decreasing function, with two different slopes, where the initial slope is bigger than that towards the end}.

$\bullet$ It is reasonable to assume that the expected number of forwards is proportional to the attractiveness factor $\rho > 0$. Hence, we model the TeF function $m(\cdot)$ as $m(a) = \rho m_N(a) ~ \forall ~ \ax$, where $m_N(\cdot)$ is the TeF with $\rho = 1$. Such a factorization of $m(\cdot)$  is supported by the experiments  on the SNAP dataset 
for a wide range\footnote{\label{footnote}Such a common  fit is good mostly for $\rho \geq 0.4$, for others we directly derive a good linear fit of $m(\cdot)$ by trial-and-error, see section \ref{sec_numerical}.} of $\rho$.
The latter function $m_N(\cdot)$ is determined solely by the characteristics of the network. Basically, $m_N(\cdot)$ corresponds to a hypothetical situation where every recipient shares with all of its friends, but new shares are only the effective forwards. Conclusively:
\begin{align}\label{eqn_network_mean_func}
m(a) &= \rho m_N(a), \mbox{ where with } \tilde{m} := \bm - \bax (\kappa_1-\kappa_2), \nonumber \\
m_N(a) &:= 
 (\bar{m} - \kappa_1 \ax)1_{ \ax \leq \bax} + (\tilde{m} - \kappa_2 \ax)1_{ \ax > \bax}.
\end{align}
Thus, the TeF is determined by four parameters, $\bm, \kappa_1, \kappa_2$ and $\bax$. Here, $\rho\bm = \rho E[\mathcal{F}(\ax_0)]$ is the expected number of forwards in the beginning, and $ \rho\kappa_1, \rho\kappa_2$ reflect the slopes of the TeF.  Observe that the slope reduces to $\rho\kappa_2 < \rho\kappa_1$ when total shares reach a particular value, $\bax$. We assume\footnote{Then, there is a possibility of exponential growth leading to virality.} $\rho \bm > 1$.

$\bullet$ We corroborate that the above model well captures the content propagation using an instance of the SNAP dataset: (a) we first obtain the TeF curve by estimating $m_N(a)$ (i.e., with $\rho = 1$) for some values of $a$;  (b) then fit a two slope curve using a naive approach; and (c) use the estimated (four) parameters to derive the solution of the ODE \eqref{eqn_ODE} and other theoretical conclusions, for posts with different $\rho$ values.

Estimating the parameters mentioned above in step (b) is crucial in any such study. For now, we estimate the same using a simple trial-and-error method, as explained in footnote \ref{footnote} and section \ref{sec_numerical}. However, more sophisticated study is essential to estimate these parameters more accurately, and we leave it for the future. As said earlier, this study aims to provide a theoretical understanding of the CP process.

\section{ODE  analysis}\label{ode_based_analysis}
Using the TeF  \eqref{eqn_network_mean_func}, we solve the ODE \eqref{eqn_ODE}. By Theorem \ref{thrm1}, this solution (a.s.) approximates the CP process. We also derive the approximate trajectories for $\ax(t)$ and $\cx(t)$, which show the desirable saturation effect. Let $\tau_s := \inf\{t: \ax(t) = \bax \}$ and $\tau_e := \inf\{t: \cx(t) = 0 \}$, the extinction time. 

\textit{In extinction sample paths (we have, $\ups_{n} = 0$),  from ODE \eqref{eqn_ODE}, $\hat{\ups}^{n}(t) = 0$ $\forall$ $t$. We now consider viral sample paths.}

\vspace{0.5mm}
\noindent \textbf{Total and current fractions:} The fractions $(\pa, \pc)$ are given by the solution of the ODE  \eqref{eqn_ODE}, which can be directly obtained using the integrating factor (IF) approach. For $t \leq \tau_e$, the closed-form (extended) solution for $\pa$ is given by: 

\vspace{-0.4cm}
{\small
\begin{align*}
    \hspace{-1mm} \pa(t) = e^{-\int_{u_1}^t (u_4 \eta(s') + 1) ds'} \hspace{-1mm} \left(\hspace{-0.8 mm}u_2 + u_3\int_{u_1}^t \hspace{-1mm} e^{\int_{u_1}^s (u_4 \eta(s') + 1) ds'} \hspace{-1mm} ds\right), \hspace{-1mm}
\end{align*}}where  constants {\small$u : = (u_1, u_2, u_3, u_4)$} in the two phases equal:

\vspace{-0.1cm}
{\small 
\begin{equation*}
u = 
\begin{cases}
\left(0, c_0, \bm \rho, \kappa_1 \rho \right), & 0 \leq t \leq \tau_s,\\
\left( \tau_s,  \pa(\tau_s), \Tilde{m} \rho , \kappa_2 \rho  \right), & \tau_s < t \leq \tau_e.
\end{cases}
\end{equation*}}

After extinction, i.e, for $t > \tau_e$, from the ODE \eqref{eqn_ODE},  $\pa(t) = \pa(\tau_e)$. Similarly, the current fraction of shares till time $t$ evaluates to the following, again using IF approach:

\vspace{-0.2cm}
{\small
\begin{eqnarray*}
    \pc(t) \hspace{-2mm} &=& \hspace{-2mm} e^{-(t-v_1)}\left(v_2 + \int_{v_1}^t  e^{(s - {v_1})}(v_3 - v_4 a(s) - 1)  ds \right), \mbox{ with} \nonumber \\
&&  \hspace{-14mm}   v  =  (v_1, v_2, v_3, v_4)  =   
    \begin{cases}
    \left(0, c_0, \bm\rho, \kappa_1\rho \right), & 0 \leq t \leq \tau_s,\\
    \left( \tau_s,  \pc(\tau_s), \Tilde{m}\rho, \kappa_2\rho \right), \hspace{-4mm}& \tau_s < t \leq \tau_e.
    \end{cases}
\end{eqnarray*}}Further,  $\pc(t) = 0$ for  all $t \geq \tau_e$.

\noindent \textbf{Trajectories of shares:} It is more relevant to analyse the trajectories corresponding to the current and total shares; and we provide approximate expressions for the same.  Recall from \eqref{eqn_ODE}, $\ax(t) = \pa(t) \eta(t)$ and $\cx(t) = \pc(t) \eta(t)$. Thus, we begin with an approximation for  $\eta(t)$, defined in \eqref{eqn_eps_n}:
\begin{align}\label{eqn_eta}
\begin{aligned}
    \eta(t) &\approx \max\{n : \gamma + \ln(n) \leq t \}= \floor{e^{t-\gamma}} \approx e^{t-\gamma},
\end{aligned}
\end{align}where $\gamma$ is Euler-Mascheroni constant. \textit{Henceforth, we use this approximation of $\eta(t)$ in all the computations.} With this approximation and using \eqref{eqn_ODE}, the ODE for $\ax(\cdot)$ is given by:
\begin{align}\label{approx_total_ode}
    \dot{\ax} &=  
     \dot{\pa}e^{t-\gamma} + \pa e^{t-\gamma}1_{\cx > 0} = m(\ax)e^{t - \gamma}1_{\cx > 0}.
\end{align}
The solution of the above ODE with $m(\cdot)$ as in \eqref{eqn_network_mean_func}, can be obtained using the standard techniques in ODE theory:
\begin{equation}\label{eqn_total_ode}
\boxed{
a(t) = 
\begin{cases}
 w_1 - w_2 e^{-w_3 e^{t}}, &\mbox{ when } 0 \leq t \leq \tau_e,\\
a(\tau_e), & \mbox{ when } t > \tau_e,
\end{cases}}
\end{equation}
where the constants are given by the following (with $w = (w_1, w_2, w_3)$, observe by continuity $a(\tau_s) = \bax$):

\vspace{-0.3cm}
{\small 
\begin{equation}\label{total_pop_constants}
w = 
\begin{cases}
\left( \frac{\bm}{\kappa_1}, (w_1 - \ax_0)e^{w_3}, \kappa_1 \rho e^{-\gamma} \right), & 0 \leq t \leq \tau_s,\\
\left( \frac{\tilde{m}}{\kappa_2}, (w_1- \bax)e^{w_3  e^{\tau_s}},  \kappa_2 \rho e^{-\gamma} \right), & \tau_s < t \leq \tau_e.
\end{cases}
\end{equation}}Observe that $\ax(\cdot)$ saturates as the current shares get extinct. We believe, only STP-BP can capture these effects. Proceeding as in \eqref{approx_total_ode}, the ODE for $\cx(t)$ is given by:

\vspace{-6mm}
\begin{align}\label{approx_current_ode}
    \dot{\cx} = (m(\ax) - 1)e^{t-\gamma}1_{\cx > 0}.
\end{align}
By solving, the trajectory of the current shares is given by:
\begin{empheq}[box=\widefbox]{align}\label{eqn_current_ode}
    \cx(t) &= \left(\cx(\varphi) - \ax(\varphi) + \ax(t) + e^{-\gamma}(e^\varphi - e^t)\right)1_{t < \tau_e}, 
\end{empheq}
where $\varphi := \tau_s 1_{t > \tau_s}$. After extinction (for $t \geq \tau_e$), $\cx(t) = 0$. 

\vspace{0.5mm}
\noindent \textbf{Shares at transition epochs:} From Theorem \ref{thrm1}, the value of the embedded chain at the transition/wake-up epoch $n$ is approximated by the ODE solution at $t_n$, defined in \eqref{eqn_eps_n}. Thus, it is more important to evaluate $\ax(t)$ and $\cx(t)$ at time points $ t = t_n$. Towards this, define $n_s := \eta(\tau_s)$ and $n_e := \eta(\tau_e)$ as the respective counterparts of $\tau_s$ and $\tau_e$. Using the same approximation ($t_n \approx \gamma + ln(n)$) as in \eqref{eqn_eta}, now for the mapping $n \mapsto t_n$, the shares, by \eqref{evolve_SA}, \eqref{eqn_total_ode} and \eqref{total_pop_constants}, equal:

\vspace{2mm}
\fbox{
\begin{minipage}{0.45\textwidth}
\hspace{-2mm}\vspace{-2mm}
\begin{align}\label{eqn_total_epoch}
\hspace{-5mm}\ax(t_n) &=
     w_1 - w_2 e^{-n w_3 e^\gamma},  \mbox{ when } 0 \leq n \leq n_e,\\
\hspace{-5mm}     \cx(t_n) &=  w_1 - w_2 e^{-n w_3 e^\gamma} - n, \mbox{ when } 0 \leq n \leq n_e.\label{eqn_current_epoch}
\end{align}
\end{minipage}}

\vspace{2mm}
\noindent In the above, \eqref{eqn_current_epoch} follows because  $\ax(t_n)-\cx(t_n) = n$ from \eqref{evolve_SA}. After extinction, $(\ax(t_n), \cx(t_n)) = (\ax(\tau_e), 0)$  $\forall$ $n \geq n_e$. 

We now investigate other important metrics related to CP.

\noindent \textbf{Growth rates:} In standard BPs, the population exhibits dichotomy, i.e., either the population grows exponentially large (at a constant rate),  or declines to zero. In the former case, both $\Cx(t_n)$ and $\Ax(t_n)$ grow like $e^{(\rho \bm -1)\tilde{\tau}_n}$, where $\tilde{\tau}_n := \lambda \tau_n$ is independent\footnote{Using properties of exponential distribution, this can easily be proved. Such an effect is seen since we are discussing the embedded chain.} of $\lambda$, as time progresses. While in STP-BP, from \eqref{eqn_total_epoch}, which again approximates $\Ax(t_n)$ a.s., and time asymptotically, it is clear that the total shares have exponential growth, however, the rates are different in the two phases. The growth rate in the initial phase is $w_3 e^\gamma = \kappa_1 \rho$, while it decreases to $\kappa_2 \rho$ in the later phase (since $\kappa_1 > \kappa_2$). Further, \textit{the current shares also experience an initial exponential growth, which is further modulated by the growing factor of $n$, in \eqref{eqn_current_epoch}, leading to an eventual linear fall}. This illustrates the modified dichotomy discussed after Theorem \ref{thrm1}. More attractive posts have higher growth rate.

At this point, we would like to admit that discussing the growth rates at transition epochs is a non-standard practice in the BP literature. However, for the saturated BP, that does not grow forever, we believe such a discussion is relevant and important. In future, we plan to include the influence of $\{\tau_n\}$ on these growth rates, by extending the ODE analysis to fractions $\left\{\nicefrac{\tau_n}{n}\right\}$; this would help us derive the standard growth patterns discussed in the BP literature for STP-BP. 
 
\noindent  \textbf{Peak of current shares:} Define $\cx^* = \sup_t \cx(t)$, i.e., the peak (maximum) current shares. It can be obtained from \eqref{eqn_current_ode} $\left(c''(t) < 0\right)$ and equals (recall $n_s = a(\varphi) - c(\varphi)$):
\begin{eqnarray}\label{peak}
    \cx^* = w_1 - \frac{\left(1 + \ln(w_2 w_3 e^\gamma) \right)}{w_3 e^\gamma}, \mbox{ as } e^{-\gamma  + \varphi} = n_s,
\end{eqnarray}
where $\varphi$ takes different values in two phases as in \eqref{eqn_current_ode}.

\noindent \textbf{Life span and Max reach:} Recall $n_e$ is the epoch at which CP of the post terminates, or in other words, it represents the life span of the post. Now, substituting $c(t_n)$ from \eqref{eqn_current_epoch}, $n_e$ can be written as the solution of the equation:

\vspace{-6mm}
\begin{eqnarray}\label{eqn_max_share}
    w_1 - w_2 e^{-n_e w_3 e^\gamma} = n_e,
\end{eqnarray}
where  $w$ is given by the second line of \eqref{total_pop_constants}. 
Further, by \eqref{eqn_current_epoch}, the \textit{max reach} or saturated total shares $\ax(t_{n_e})$ equal $n_e$.

\noindent \textbf{Probability of virality:} When one transitions from super-critical to sub-critical regime, it is clear that the process gets extinct a.s. (see modified dichotomy). One of the important questions related to CP is the chances of virality, where total/current shares grow significantly large. This necessitates the definition of a different probability related to STP-BP, the probability of virality (PoV). In view of the initial exponential growth of current shares \eqref{eqn_current_epoch} and the dichotomy remark after Theorem \ref{thrm1}, we define PoV, $p_\Delta := P(\Cx(t) > \Delta \mbox{ for some } t)$ for some threshold $\Delta > 0$. We conjecture that for small $\Delta$, $p_\Delta \approx 1 - p_e$, where $p_e$, the probability of extinction of standard population independent BP solves the equation: $f(s) = s$, where $f(\cdot)$ is the PGF of $\Gamma(a_0)$ (see \cite{athreya2001branching}). This is an important aspect for future study.

\begin{figure*}
\vspace{2mm}
\begin{minipage}{0.33\textwidth}
    \includegraphics[width = 6.1cm, height = 4.2cm]{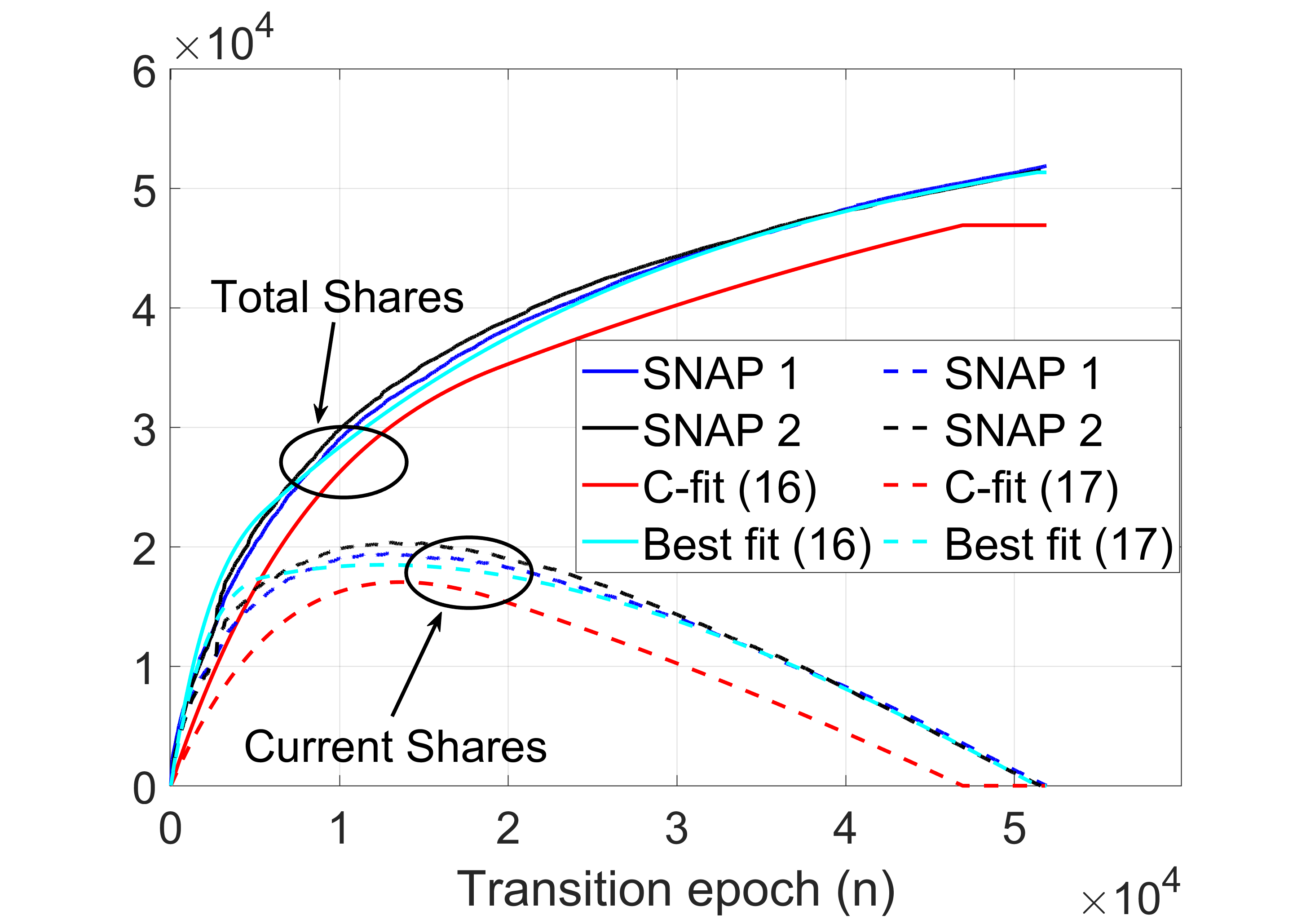}
\end{minipage}%
\begin{minipage}{0.33\textwidth}
    \includegraphics[width = 6.1cm, height = 4.2cm]{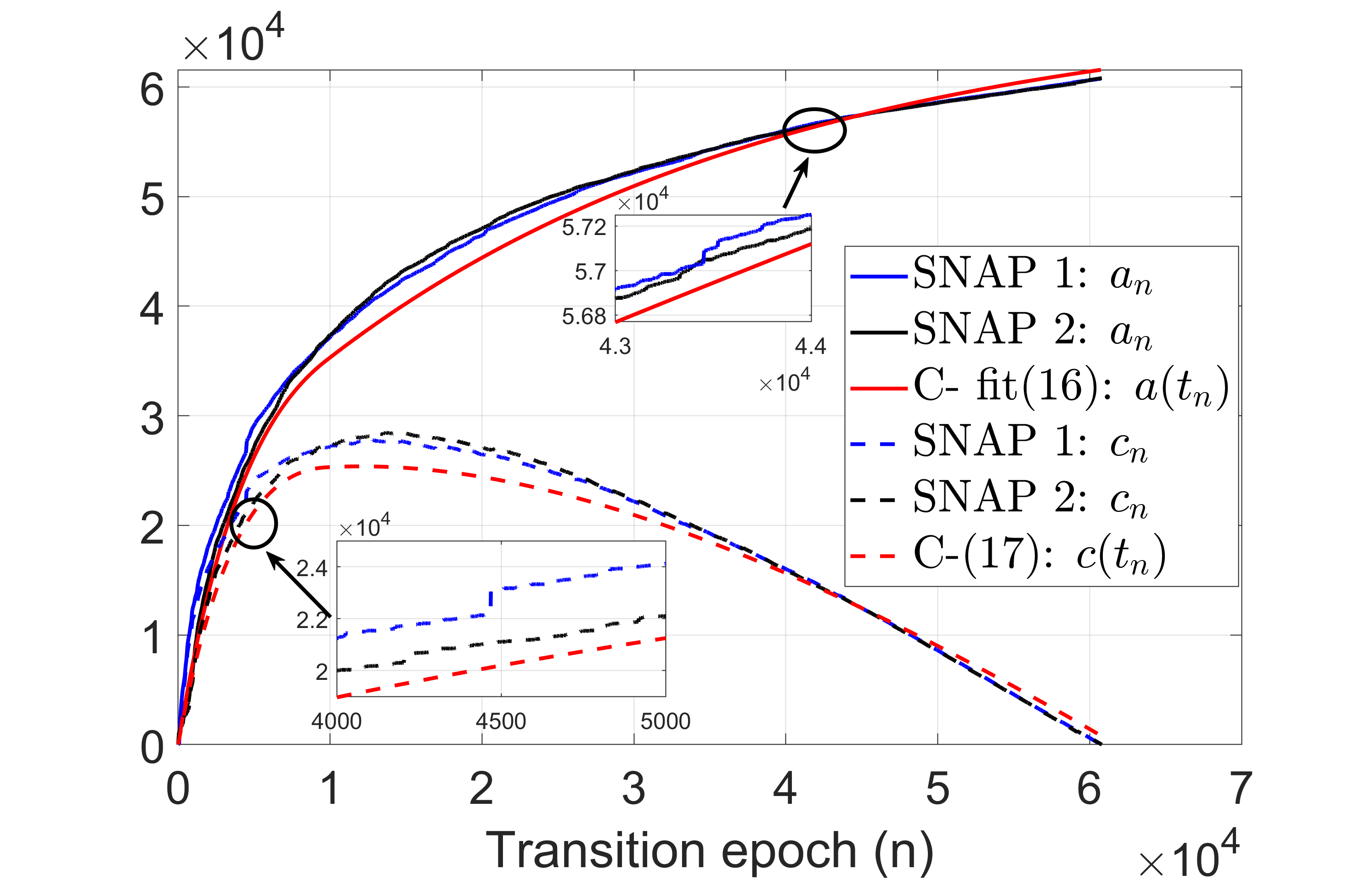}
\end{minipage}%
\begin{minipage}{0.33\textwidth}
    \includegraphics[scale = 0.27]{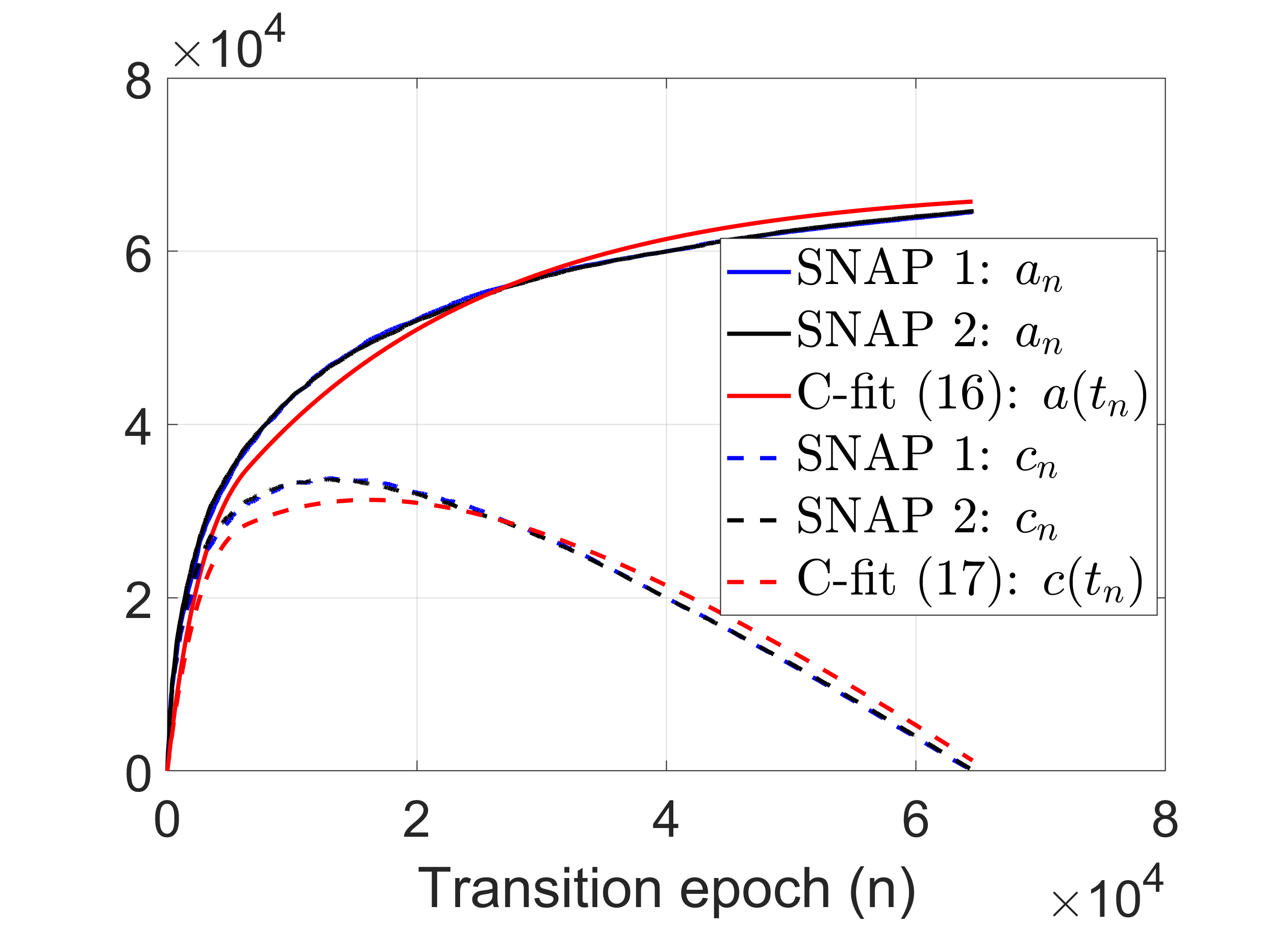}
\end{minipage}
\caption{\label{fig_SNAP}Validation of theoretical trajectories against two instances of CP on SNAP dataset, for $\rho = 0.2, 0.4, 0.6$ respectively}
\vspace{-4mm}
\end{figure*}

\section{Numerical Experiments}\label{sec_numerical}
We now perform exhaustive Monte-Carlo (MC) simulations on the SNAP Twitter dataset \cite{mcauley2012learning} to validate our theory. The dataset consists of inter-connections among $81,306$ users, with $29.77$ average number of friends. \textit{In all the case studies, we consider only viral sample paths.}

\begin{figure}[htbp]
\vspace{-3mm}
\begin{minipage}{0.5\linewidth}
    \includegraphics[scale = 0.2]{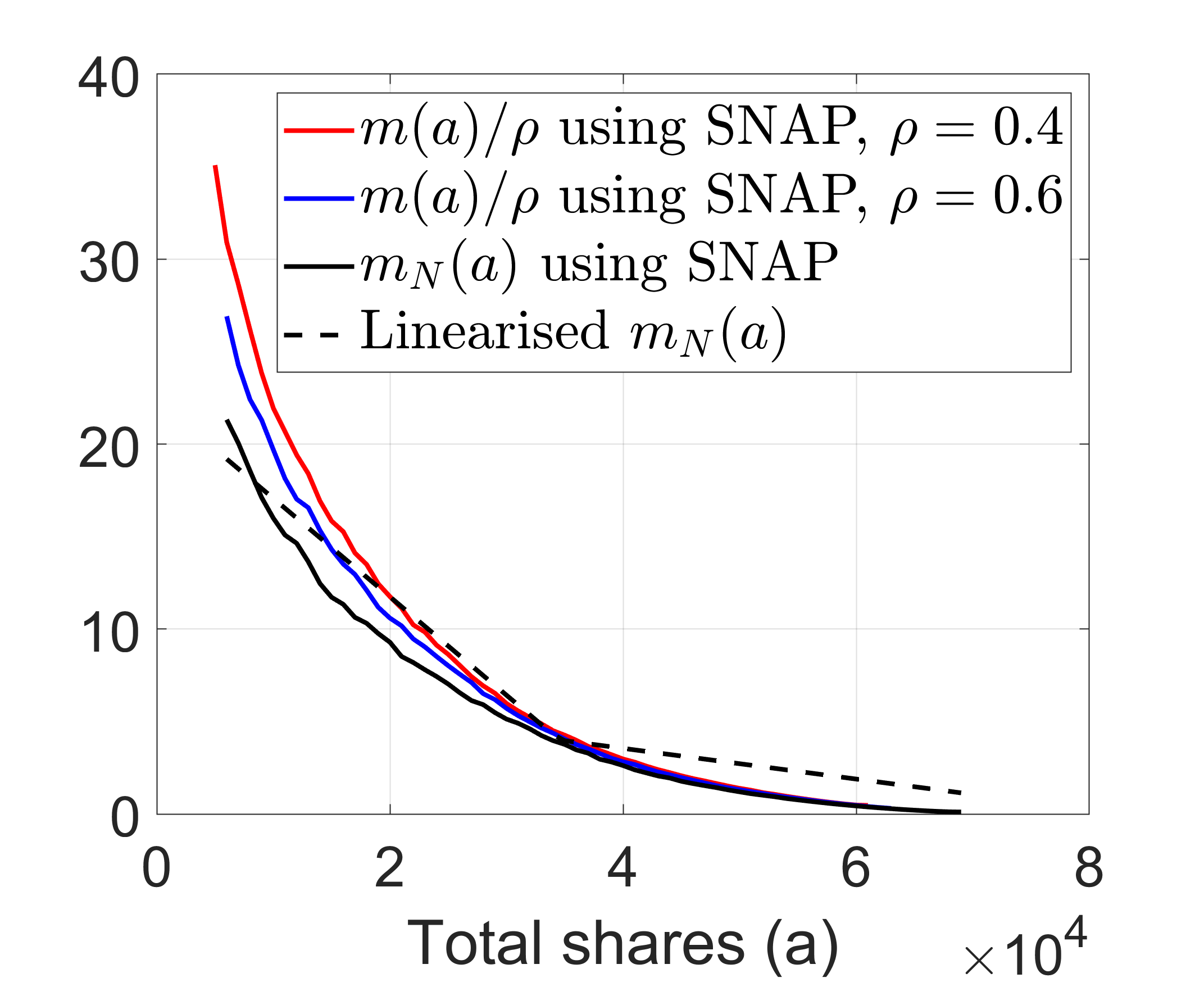}
\end{minipage}%
\begin{minipage}{0.5\linewidth}
\includegraphics[scale = 0.18]{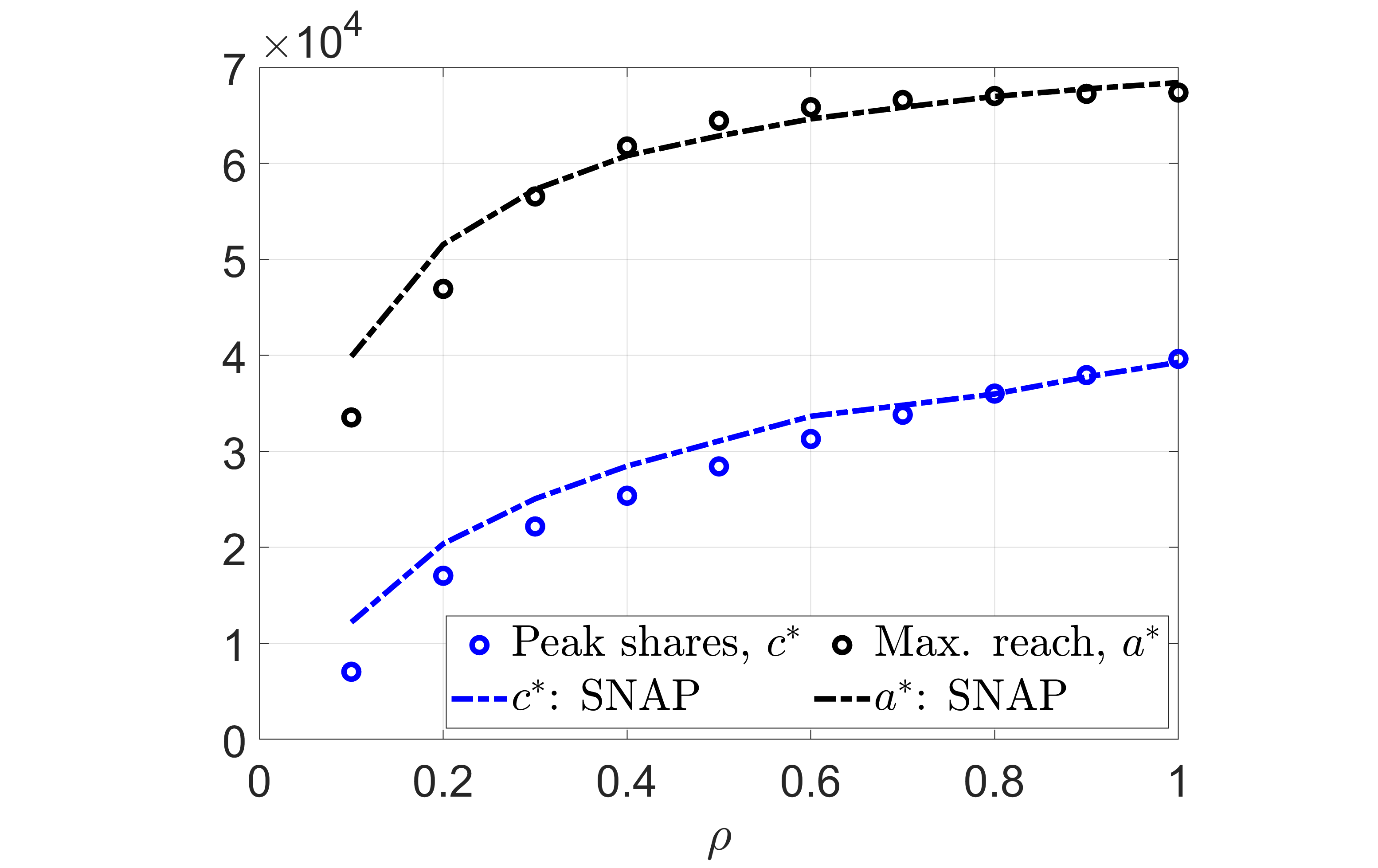}
\end{minipage}
\caption{(i) Left: piece-wise linear TeF; (ii) Right: Peak shares \eqref{peak} and Max reach  \eqref{eqn_max_share} versus SNAP estimates \label{fig_synthetic_and_mean_curve}}
\vspace{-1mm}
\end{figure}

\hide{
\noindent \textbf{Validation using synthetic data:} In Theorem \ref{thrm1}, we established that the solutions of the ODE \eqref{eqn_ODE} can almost surely approximate the dynamics of the CP process over any finite time window. We corroborate this result using MC simulations for  $\bax = 1400$, $\bm \rho = 5$,  $\rho \kappa_1 = 0.003$, and $\rho \kappa_2 = 28 \times 10^{-5}$ in Fig. \ref{fig_synthetic_and_mean_curve}(i). In particular, we plot two sample paths of the embedded chain \eqref{evolve_SA}, against the Piccard's iterates of the ODE \eqref{eqn_ODE}, approximate trajectories derived in \eqref{eqn_total_ode}, \eqref{eqn_current_ode} and the trajectories obtained at transition epochs in \eqref{eqn_total_epoch}, \eqref{eqn_current_epoch} for the total and current shares. In the said figure, we apply Piccard's iterates with $n_m = 5$. Then, we evaluate the shares (represented by red triangles) by multiplying the fractions obtained by $\eta(t_k)$ for finite $k > n_m$. Similarly, the approximate trajectories are plotted at transition epochs, and are displayed by red circles. 
From the figure, it can be seen that Theorem \ref{thrm1} indeed holds true for synthetic data, and further the derived trajectories also well approximate the embedded chain. 
}

\noindent \textbf{MC simulations over the dataset:} Say the content provider initially shares the post with $2$ random seed users chosen from the dataset, identified by their user IDs. These users are added to the total and current-shares lists. At any time, one random user from the current-shares list wakes-up and forwards the post to a random subset of its friends, each chosen independently with probability $\rho$. Out of this subset, we ignore the friends who had the post. Then, we delete this user from the current-shares list, and update the two lists with the new effective forwards. The propagation continues in this manner and terminates when the current-shares list is empty, thus completing one sample run of the CP process.

\noindent \textbf{Estimation of TeF:} 
To estimate the TeF discussed in section \ref{sec_mean}, we create bins of equal length ($1000$ in our case). Then, we simulate the CP process $861$ times over the dataset  for $\rho = 1$. In each run and bin $n$, we count: (i) the number of transitions that have occurred and (ii) the sum total of effective forwards,  while the number of total shares belongs to that bin. These two entries are accumulated bin-wise over $861$ viral sample paths. For each bin, we divide the effective forwards by the number of transitions, which  represents the estimated TeF curve (see solid black curve in Fig. \ref{fig_synthetic_and_mean_curve}(i)).

We repeat the above routine for $\rho = 0.4, 0.6$ as well, and the corresponding estimates of $m(\cdot)$ are plotted in Fig. \ref{fig_synthetic_and_mean_curve}(i) after dividing by the respective $\rho$ values. The resultant picture  gives the confidence that one can derive the individual $m(\cdot)$ curves by using $m_N(\cdot)$ as suggested by \eqref{eqn_network_mean_func}. We also plot an approximate piece-wise linear curve (see the dashed curve in Fig. \ref{fig_synthetic_and_mean_curve}(i)), with the parameters $\bm  = 21.321042, \kappa_1 = 532\times 10^{-6}, \kappa_2 = 83\times 10^{-6}$ and $\bax = 35000$ obtained using trial and error method. Henceforth, we refer this $m_N(\cdot)$ curve as the \textit{common-fit (C-fit) curve}. For some sets of simulations, again using the trial and error method, we individually choose \textit{best fit $m(\cdot)$ curve} for the given case study. 

In Fig. \ref{fig_SNAP}, we obtain the theoretical curves for total and current shares using C-fit TeF (discussed above), which are then compared with two instances of the CP process on the SNAP dataset. The theoretical curves well approximate the dataset curves for $\rho = 0.4, 0.6$. In fact, we observed this is true for $\rho \geq 0.4$, in general. However, for lower values of $\rho$ (e.g., $\rho = 0.2$ in Fig. \ref{fig_SNAP}), the best fit curve better approximates. This suggests that one can obtain C-fit curves for higher and lower ranges of $\rho$ separately. At last in Fig. \ref{fig_synthetic_and_mean_curve}(ii), for different $\rho$ values, we show that the peak shares and max. reach from the C-fit curve approximate the respective values for an instance of CP on the dataset with a maximum error of $11.4645\%$ and $2.4837\%$ respectively (for $\rho \geq 0.4$). 

\section{Conclusions}
In this paper, we studied the saturation effect experienced by the total shares/copies of the post due to continual re-forwarding on OSNs. We captured the dynamics via newly introduced, saturated total-population dependent branching process. The analysis uses the stochastic approximation technique, which provides an ODE dependent on the post's total expected (effective) forwards, TeF. We modelled TeF as a piece-wise linearly decreasing function with two slopes. The derived trajectories (dependent on four network-specific parameters) asymptotically and almost surely approximate the total and unread copies over any finite time interval. 

Unlike classical dichotomy (either explosion or extinction), the unread copies observe explosion followed by extinction or direct extinction under saturated BPs. Interestingly, maximum reach (number of users that received the post) and the life span (number of users that read the post) are equal. We showed that derived expressions provide a good fit to the simulated instances of propagation over the SNAP dataset. Here, we fitted the parameters of TeF using a naive method; however, advanced methods of estimation/learning can improve the approximation.

\section{Appendix}

\begin{lemma}\label{lemma_existence}
There exists a unique, extended, continuous solution of the ODE \eqref{eqn_ODE} over any finite time interval. 
\end{lemma}
\begin{proof}
\arxiv{
The proof follows using existence results for Carath\'{e}odory conditions given in \cite[Theorem 1.1, 1.3, pp.  43, 47]{coddington1955theory},  \cite[Lemma 2.1]{singer2006bounding} and the uniqueness follows by Lemma \ref{lemma_unique_ode_sol}. The rest of the proof is completed in \cite{arxiv}.
}{At first, we show that if at all there exists a solution for the non-smooth, non-autonomous ODE \eqref{eqn_ODE}, then, the ODE solution is bounded. Towards this, say ${\psi}^c(0) = {\psi}^a(0) = a_0$. If possible, say there exists a solution $\hat{\ups}(\cdot)$ for the ODE \eqref{eqn_ODE} for all $t < T$ (for some fixed, finite $T$). Consider the ODE for $z(\cdot) := (z^c(\cdot), z^a(\cdot))$ (see assumption \textbf{(A)}):
\begin{align}\label{eqn_z_ode}
    \dot{z}^c = b_1 - 1,\mbox{ and } \dot{z}^a = b_1, \mbox{ with } E[\hat{\Gamma}_1] =: b_1.
\end{align}The RHS of the above ODE is smooth and autonomous.
Further, consider $\bar{t}$ such that $\eta(t) = n$ for all $t \leq \bar{t}$. Then, for $t \leq \bar{t}$, we have:
\begin{align}
    \dot{\pc} &= m(\pa n) - 1 - \pc, \mbox{ and } \dot{\pa} = m(\pa n) - \pa.
\end{align}It is clear that the RHS of the above ODE is Lipschitz continuous and autonomous, and further that $\dot{\ups } < \dot{z}$,  for all $t \leq \bar{t}$. Hence, if $z(0) = \ups(0)$, then, $\hat{\ups}(t) < z(t)$ for all $t \leq \bar{t}$ (e.g., see \cite[pp. 168]{piccinini2012ordinary}).

Thus,  \cite[Lemma 2.1(ii)]{singer2006bounding} is not true (while remaining hypotheses are true), and hence by \cite[Lemma 2.1(i)]{singer2006bounding}, we have the following:
\begin{align}\label{eqn_bound_ode}
    \hat{\psi}^a(t) &< z^a(t) 
    < a_0 + b_1 T =: \beta, \mbox{ and }\\
    \hat{\psi}^c(t) &< z^c(t) 
    < \beta, \mbox{ for all } t \in [0, T]. \nonumber
\end{align}
Next, consider the set $\cd := [0, T] \times [-2\beta, 2\beta]^2$. Then, under assumption \textbf{(A)} and using \cite[Theorem 1.3, pp.  47]{coddington1955theory}, there exists 
a solution $\hat{\ups}(\cdot)$ for all $
t < \tau$, where 

\vspace{-4mm}
{\small
$$
\tau:= \min\bigg\{\inf_t\{\hat{\psi}^c(t) \in \{-2\beta, 2\beta\}, \hat{\psi}^a(t) \in \{-2\beta, 2\beta\}\}, T\bigg\}.
$$}
Lastly, note that $\tau = T$; this is so because $\tau < T$ contradicts \eqref{eqn_bound_ode} as $\ups \geq 0$.\footnote{One can prove that $\ups \geq 0$ using \eqref{eqn_ODE} and simple lower-bounding arguments like above.} This proves the existence of the solution of the ODE for all $t \in [0, T]$. 

Next, the uniqueness of the solution holds by Lemma \ref{lemma_unique_ode_sol}. Lastly, the continuity for the solution follows by the integral representation of the solution, and because the RHS of the ODE \eqref{eqn_ODE} can be bounded by $b_1 + \beta$. }
\end{proof}

\noindent \textbf{Proof of Theorem \ref{thrm1}:} \textbf{Part (i)} The proof of this part follows closely as in \cite[Theorem 2.1, pp. 127]{kushner2003stochastic}, but the RHS of the ODE in our case  is only measurable. 
Let $n \geq 0$. Using \eqref{eqn_SA}, one can re-write $\Ups_n$ as:

\arxiv{\vspace{-5mm}}{}
\begin{align}\label{eqn_scheme2}
\begin{aligned}
 \Ups_{n+1} &= \Ups_n + \epsilon_n L_n, \mbox{ where } L_n := ( L_n^{ c}, L_n^{ a}), \mbox{ and}\\
 L_n^{c} &:=  \left[\offs_{n}(A_{n-1}) - 1 - \Pc_{n-1} \right ]1_{\Pc_{n-1} > 0}, \\
L_n^{ a} &:=  \left [\offs_{ n}(A_{n-1})  - \Pa_{n-1} \right ]1_{\Pc_{n-1} > 0}.
\end{aligned}
 \end{align}
 
\textbf{Interpolated trajectory.} Let $\Ups^n(\cdot) = (\Psi^{n, c}(\cdot), \Psi^{n, a}(\cdot))$ be the piece-wise interpolated trajectory defined as (see \eqref{eqn_scheme2}):

\arxiv{\vspace{-5mm}}{\vspace{-4mm}}
{\small
\begin{align}\label{eqn_piecewise}
    \Ups^{n}(t) &= \Ups_n + \sum_{i=n+1}^{\eta(t_n +t)}(\Ups_i - \Ups_{i-1}) \nonumber \\
    &= \Ups_n + \sum_{i=n}^{\eta(t_n +t)-1}\epsilon_{i} L_{i}, \mbox{ for any }t \geq 0.
\end{align}}Let  $\bar{g}(\Ups_n, n) := E[ L_n | \mathcal{F}_n]$, i.e., the conditional expectation of $L_n$ with respect to ${\cal F}_n  := \sigma\{\Ups_k : 1 \leq k < n \}$ and $\delta M_n :=  L_n - \bar{g}(\Ups_n, n)$. Then, \eqref{eqn_piecewise} can be re-written component-wise as (for each $k \in \{a, c\}$):

\vspace{-4mm}
{\small
\begin{align}\label{eqn_linear2}
    \Psi^{n, k}(t) &= \Psi^k_n + \sum_{i=n}^{\eta(t_n +t)-1}\epsilon_{i} \left(\delta M_{i}^k + \bar{g}^k(\Ups_i, i)\right) \nonumber \\
    &\hspace{-4mm}= \Psi^k_n + \int_0^t \bar{g}^k(\Ups^n(s), n) ds + \varepsilon^{n, k}(t), \mbox{ where}\\
    \varepsilon^{n, k}(t) &:= M^{n, k}(t) + \rho^{n, k}(t) \mbox{ with }
    M^{n, k}(t) := \sum_{i=n}^{\eta(t_n +t)-1}\epsilon_i \delta M_i^k, \nonumber  \\ 
    \rho^{n, k}(t) &:=  \sum_{i=n}^{\eta(t_n +t)-1} \hspace{-2mm}\epsilon_{i} \bar{g}^k(\Ups_i, i) - \int_0^t \bar{g}^k(\Ups^n(s), n) ds. \nonumber
\end{align}}It is important to note  that $\bar{g}(\Ups_n, n)$ is the RHS of the ODE \eqref{eqn_ODE} (as $\eta(t_n) = n$). 

Next, we begin by proving that the BP trajectory (see \eqref{eqn_SA}) can be bounded (under assumption \textbf{(A)}) as follows:

\vspace{-4mm}
{\small
\begin{align*}
0 \leq \Psi^{n, c}(0) = \Psi_n^c &\leq \frac{1}{n} \left(\sum_{k = 1}^n \Gamma_k(A_{k-1})1_{\Pc_{k-1} > 0} + \ax_0 \right)\\
&\leq \frac{1}{n} \left(\sum_{k = 1}^n \hat{\Gamma}_k + \ax_0 \right) := \hat{\Pi}_n.
\end{align*}}By strong law of large numbers,  $\hat{\Pi}_n \to E[\hat{\Gamma}_1]$ a.s.
Consider any  such sample path ($\omega$). Then, for any $\epsilon > 0$, there exists $N_\epsilon(\omega)$ such that:
\begin{align}\label{eqn_bound}
\Psi^{n, c}(0) &\le \hat{\Pi}_n \leq M(\omega) \mbox{ for all } n, \mbox{ where } b_1 := E[\hat{\Gamma}_1], \mbox{ and} \nonumber \\
M(\omega) &:= \max\{\max\{\hat{\Pi}_i : 0 \leq i < N_\epsilon(\omega)\}, b_1 + \epsilon\}.
\end{align}

Now, we will prove that $M^n(\cdot) = (M^{n, c}(\cdot), M^{n, a}(\cdot))$ and $\rho^n(\cdot)$ individually converge to $0$ (as $n \to \infty$) uniformly on any bounded interval\arxiv{ over an almost sure subset, $N$ of $\{\hat{\Pi}_n \to b_1\}$ (i.e., with $P(N) = 1$); these steps follow majorly as in \cite[Theorem 2.1, pp. 127]{kushner2003stochastic} and are provided in \cite{arxiv}. Thus, we have the convergence of $\varepsilon^n(\cdot)$ under the norm \eqref{eqn_norm}, completing part (i). }{. It suffices to prove uniform convergence for sample paths $\omega \in \{\hat{\Pi}_n \to b_1\}$. 
We prove the claim for $\pc$-component, and it can proved analogously for the $\pa$-component as well. Henceforth, the convergence will be proved w.r.t. $n$, where ever not mentioned explicitly. 
 
Now, define $M_n^{c} := \sum_{i=0}^{n-1}\epsilon_i \delta M_i^{c}$. Then, it is easy to prove that $(M_n^{c})$ is a Martingale   with respect to $(\mathcal{F}_n)$. Thus, using Martingale inequality, 
for each $\mu > 0$ (as in \cite[Theorem 2.1, pp. 127]{kushner2003stochastic}), with $E_n(\cdot)$ denoting the expectation conditioned on $(\mathcal{F}_n)$: 

\vspace{-4mm}
{\small
$$
P\left\{\sup_{m\leq j \leq n} |M_j^{ c} - M_m^{ c}| \geq \mu \right\} \leq \frac{E_n\left[\left(\sum_{i=m}^{n-1} \epsilon_i \delta M_i^{ c} \right)^2 \right]}{\mu^2}.
$$}
Observe, $E\left[\delta M_i^{ c} \delta M_j^{ c}\right]  = 0$ for $i < j$. Let $O(\omega)$ be the upper-bound on the ODE solution for $t \in [0, T]$, see Lemma \ref{lemma_existence}. Then, from \eqref{eqn_ODE} and \eqref{eqn_bound}, 
\begin{align}\label{eqn_g_bound}
    |\bar{g}^i(\ups(\cdot), \cdot)| < b_1 + 1 + O(\omega) \mbox{ for each } i \in \{a, c\}.
\end{align}Thus,  under \textbf{(A)}, $\sup_n E_n|L_n^{ c}- \bar{g}^c(\Ups_i, t_i)|^2 < K$ for some finite $K$.
Using this,  we have: 

\vspace{-4mm}
{\small
\begin{align*}
    P\left\{\sup_{m\leq j \leq n} |M_j^{c} - M_m^{ c}| \geq \mu \right\} &\leq \frac{\sum_{i=m}^{n-1} \epsilon_i^2 E_n\left| \delta M_i^{ c} \right|^2}{\mu^2} \nonumber \\
    &\hspace{-30mm}= \frac{\sum_{i=m}^{n-1} \epsilon_i^2 E_n\left| L_i^{ c} - \bar{g}^c(\Ups_i, i) \right|^2}{\mu^2} \leq \frac{K}{\mu^2} \sum_{i=m}^{\infty} \epsilon_i^2.
\end{align*}}
By first letting $n \to \infty$ (and using continuity of probability), then, letting $m \to \infty$, for each $\mu > 0$, we have:
\begin{align}
    \lim_{m \to \infty} P\left\{\sup_{m\leq j } |M_j^{ c} - M_m^{ c}| \geq \mu \right\} &= 0. \label{eqn_equi_cont_M}
\end{align}
Now, by \eqref{eqn_equi_cont_M} and continuity of probability for each $k > 0$, 
$P(A_k) = 1$, where $A_k := \lim_{m \to \infty} \sup_{m\leq j } |M_j^{c} - M_m^{ c}| < 1/k$. We further restrict our attention to  sample paths  $\omega \in  N:= (\cap_k A_k) \cap \{\hat{\Pi}_n \to b_1 \}$. For any such $\omega$, using \eqref{eqn_linear2}:
\begin{align*}
&\sup_{t\geq 0}|M^{n, c}(t)| = \sup_{t \geq 0} \left|M^{ c}_{\eta(t_n+t)} - M^{c}_n \right|  
= \sup_{j \geq n} |M^{ c}_{j} - M^{c}_n|.
\end{align*}This implies:

\vspace{-4mm}
{\small
\begin{align*}
\lim_{n \to \infty} \sup_{t \in [0, T]}|M^{n, c}(t)|  &\leq \lim_{n \to \infty} \sup_{t \in [0, T]}\left|\sum_{i=n}^{\eta(t_n+t)}\epsilon_i \delta M_i^k\right| \\
&\hspace{-20mm}\leq   \lim_{n \to \infty} \sup_{\eta(t_n+t) + 1 \geq n} |M^{ c}_{\eta(t_n+t) + 1} - M^{c}_n| < 1/k.
\end{align*}}Letting $k \to \infty$, we get, $M^{n,  c}(\cdot) \to 0$ uniformly on each bounded interval. 

For $\rho^{n, c}(\cdot)$, note that for  $t = t_k - t_n$ $(k > n)$, $\rho^{n, c}(t) = 0$. Thus, under \eqref{eqn_g_bound}, for any $|t| \leq T$ (as $\epsilon_{\eta(t_n + t)} \le \epsilon_n$):

\vspace{-4mm}
{\small
\begin{align*}
    |\rho^{n, c}(t)|
    &\leq  \int_{t_{\eta(t_n+t) } - t_n}^t \left|\bar{g}^c(\Ups^n(s), \eta_n) \right|ds 
    <    \epsilon_n (b_1+1+O).
\end{align*}}Thus, $\rho^{n, c}(\cdot) \to 0$ uniformly on each bounded interval. }

\textbf{Part (ii)}
We construct this proof using Maximum Theorem, which provides parameterized continuity of the optimizers. We begin by constructing the required elements  (i.e., appropriate objective function and domains).

\textbf{Ingredients for Maximum Theorem.} Fix any $\omega \in N$. \arxiv{Let $O(\omega)$ be the upper-bound on the ODE solution for $t \in [0, T]$, see Lemma \ref{lemma_existence}.}{} \arxiv{Thus,}{Then,} the interpolated trajectory $\Ups^n(\cdot)$ and the ODE solution $\hat{\ups}^n(\cdot)$ are bounded as (see \eqref{eqn_bound}):

\vspace{-4mm}
{\small
$$
\sup_{t} \Ups^n(t) = \sup_{n} \Ups_n < 1.1 M(\omega), \mbox{ and } \sup_{t \in [0, T]} \hat{\ups}^n(t) < 1.1 O(\omega).
$$}
With the norm \eqref{eqn_norm}, let $\cd^2$ be the Banach space of all those $\ups(\cdot)$ such that both $\pc, \pa$ are left continuous with right limits on $[0, T]$ and $||\ups||< \infty$. Further, let $\cd^2_B$ be the space of all those $\ups(\cdot) \in \cd^2$ such that $||\ups|| \leq C(\omega) := 1.1(M(\omega) + O(\omega))$. 
Define $\cd_p :=  \cd^2_B \times \mathbb{R}^2 \times \mathbb{R}$, and then, define the function $F(\ups; \varepsilon, u_0, \eta): \cd^2_B \times \cd_p \to \mathbb{R}$ as:

\vspace{-4mm}
{\small
\begin{align*}
F(\ups; \varepsilon, u_0, \eta) &:= \sum_{i \in \{a, c\}} \int_0^T \bigg(\Psi^i(t) - h^i(\ups; \varepsilon, u_0, \eta)(t) \bigg)^2 dt,
\end{align*}}where for any $t$, the function $h^i$ is defined as:

\vspace{-4mm}
{\small
\begin{align}\label{eqn_h}
    h^i(\ups; \varepsilon, u_0, \eta)(t) &:= u_{0, i} + \int_0^t \bar{g}^i(\ups(s), \eta) ds + \varepsilon^i(t).
\end{align}}
We prove the required continuity via the parametric continuity of the following optimization problem:

\vspace{-4mm}
{\small
\begin{align}\label{eqn_opt}
F^*(\varepsilon, u_0, \eta) := \inf_{\ups\in \cd^2_B} F(\ups; \varepsilon, u_0, \eta)  \ \forall \ (\varepsilon, u_0, \eta) \in \cd_p.
\end{align}}
It is clear that the minimizer ($\ups^*$) of \eqref{eqn_opt} is the fixed point of  the operator $\ups \mapsto h(\ups; \cdot, \cdot, \cdot)$, if one exists, and then, $F(\ups^*; \cdot, \cdot, \cdot) = 0$. Also, from \eqref{eqn_linear2}, $\Ups^n(\cdot)$ is the optimizer of \eqref{eqn_opt} at parameters $(\varepsilon, u_0, \eta) = (\varepsilon^n, \ups_n, n)$, by choice of $C(\omega)$ and domain $\cd^2_B$. Similarly, the ODE solution $\hat{\ups}^n(\cdot) \in \mbox{arginf}_{\ups\in \cd^2_B} F(\ups; 0, \ups_n, n)$, again by choice of $C(\omega)$ and domain $\cd^2_B$. We complete the remaining proof in two steps. 

$\bullet$ $\mathbf{F(\ups; \varepsilon, u, \eta)}$ \textbf{is jointly continuous}, i.e., 
if $||\ups^n - \ups|| \to 0$, $u_n \to u$, $\eta_n \to \eta$ and $||\varepsilon^n - \varepsilon||\to 0$, we have, $F(\ups^n; \epsilon^n, u_n, \eta_n) \to F(\ups; \epsilon, u, \eta)$. \arxiv{From \eqref{eqn_ODE} and \eqref{eqn_bound},  $\bar{g}^i(\ups(\cdot), \eta_k) \leq b_1 + 1 + O(\omega)$ for each $i\in\{a,c\}$}{Recall from \eqref{eqn_g_bound}, $\bar{g}^i(\ups(\cdot), \eta_k) \leq b_1 + 1 + O(\omega)$ for each $i\in\{a,c\}$}. Further, by assumption \textbf{(B)}, we have, $m(\ups^{a, n}(s)\eta_n) \to m(\ups^a(s)\eta)$.
This implies $\bar{g}(\ups^n(s), \eta_n) \to \bar{g}(\ups(s), \eta)$. Then, by applying bounded convergence theorem twice, we have the claim. 

$\bullet$ $\mathbf{\cd^2_B}$ \textbf{is weak-compact.}
Consider the projection, $p_s^i(\ups) := \ups^i(s)$, for each $i \in \{a, c\}$ and $s \in [0, T]$. For each $s, i$, we have,  $p_s^i(\cd_B^2) = [-C(\omega), C(\omega)]$, which are clearly compact. By Tychonoff's Theorem, 
$\cd_B^2$
is weak-compact under the well known product topology on $\cd^2$.

Thus, the parametric optimization problem in \eqref{eqn_opt} satisfies the hypothesis of Berge's maximum theorem (e.g., \cite{feinberg2014berges}). So, the set of optimizers defined by (for all  $(\varepsilon, u, \eta)\in \cd_p$):

\vspace{-4mm}
{\small\begin{align}\label{eqn_h_star}
{\cal H}^*(\varepsilon, u, \eta) &:= \mbox{arg inf}_{\ups \in \cd_B^2} F(\ups; \varepsilon, u, \eta) = \{\ups^*(\varepsilon, u, \eta)\} 
\end{align}}
is upper semi-continuous correspondence on $\cd_p$. 

Next, define the set $\Theta \subset \cd_P$ such that
\begin{align}\label{eqn_Theta}
\Theta := \{(\varepsilon^n, \Ups_n,  n), (0, \Ups_n, n) \mbox{ for all } n \}.
\end{align}
By Lemma \ref{lemma_unique_ode_sol}, the optimizers are unique when restricted to $\Theta \subset \cd_p$.  Thus,  ${\cal H}^*$ of \eqref{eqn_h_star} is  continuous on $\Theta$, when viewed as a function. In other words, when arguments (particularly, $(\varepsilon^n, \ups_n, n)$ and $(0, \ups_n, n)$) of ${\cal H}^*$ are close-by,  then the corresponding values of ${\cal H}^*$ are also close-by. Also, by part (i), these arguments of ${\cal H}^*$ are closing-in, as $n \to \infty$.
\eop

\begin{lemma}\label{lemma_unique_ode_sol}
The optimizer of the problem in \eqref{eqn_opt} is unique if $F^* = 0$. This also implies, the solution for ODE \eqref{eqn_ODE} is unique for any given initial condition over any bounded interval.
\end{lemma}
\begin{proof} If $F^*= 0$ for some parameter say $(\varepsilon, u_0, \eta)$, then by definition, any optimizer is a fixed point for $(h^a, h^c)$, see \eqref{eqn_h}. If possible, let $\ups_1$ and $\ups_2$ be two such distinct fixed points. Then, for each $i \in \{a, c\}$ and $j \in \{1, 2\}$, we have:
\begin{align*}
    \psi_j^i(t) &= u_0 + \int_0^t \bar{g}^i(\ups_j(s),  \eta) ds + \varepsilon^i(t) \mbox{ for any } t\geq 0.
\end{align*}
Let $k_m$ be the Lipschitz constant for function $m(\cdot)$ (see assumption \textbf{(B)}). Then, we have (see \eqref{eqn_ODE}):

\vspace{-4mm}
{\small
\begin{align}\label{eqn_pc}
    |\psi^{c}_1(t) - {\psi}^{c}_2(t)| &= \left|\int_0^t \bigg(\bar{g}^c(\ups_1(s), \eta) - \bar{g}^c(\ups_2(s), \eta) \bigg) ds \right| \nonumber\\
    &\hspace{-18mm}\leq \int_0^t \left|  m(\psi^{a}_1(s) \eta) - m(\psi^a_2(s) \eta)\right| ds + \int_0^t \left| \psi^{c}_1(s) - {\psi}^{c}_2(s)  \right| ds \nonumber \\
    &\hspace{-18mm}\leq  (k_m \eta +1)  \int_0^t u(s) ds.
\end{align}}Similarly,

\vspace{-10mm}
{\small
\begin{align}\label{eqn_pa}
    \hspace{6mm}|\psi^a_1(t) - {\psi}^a_2(t)| 
    \leq (k_m \eta +1) \int_0^t  u(s) ds.
\end{align}}
Define $u(s) := \max\{\left| \psi_1^i(s) - \psi_2^i(s) \right|: i \in \{a, c\}\}$ for each $s \geq 0$. Then, from \eqref{eqn_pc}, \eqref{eqn_pa}, we have:
\begin{align*}
    u(t) &\leq \sum_{i \in \{a, c\}}\left| \psi_1^{i}(s) - \psi_2^{i}(s) \right| \leq 2(k_m \eta+1) \int_0^t u(s) ds.
\end{align*}
Applying Gronwall inequality, we have $u(t) = 0$ for each $t \in [0, T]$. This implies, $||\psi_1^{a} - \psi_2^{a}|| = 0$ and $||\psi_1^{c} - \psi_2^{c}|| = 0$, i.e., $||\ups_1 - \ups_2|| = 0$.
\end{proof}

\arxiv{\bibliographystyle{IEEEtran}
\bibliography{references_camera}
}{

}

\end{document}